\documentclass{amsart}
\usepackage{bbm,amsmath,amsfonts,amssymb}

\newtheorem{theorem}{Theorem}[section]
\newtheorem{lemma}[theorem]{Lemma}

\newtheorem{proposition}[theorem]{Proposition}
\newtheorem{definition}[theorem]{Definition}

\newtheorem{remark}[theorem]{Remark}

\newcounter{figures}[section]

\def\bZ{{\mathbb Z}}

\def\bN{{\mathbb N}}

\def\cB{\mathcal{B}}

\def\cF{\mathcal{F}}

\def\cJ{\mathcal{J}}

\def\cM{\mathcal{M}}

\def\cP{\mathcal{P}}

\def\cS{\mathcal{S}}

\def\supp{{\rm{supp}\, }}

\def\R{\mathbb{R}}

\def\cs{{c_\sharp}}
\def\fN{{\cP}}
\def\kr{{s_r}}
\def\dist{{\rm{dist}}}
\def\kkl{{k_0}}
\def\kl{{\kappa_0}}
\def\kll{{\kappa_1}}
\def\oomega{\phi}
\def\KK{K}
\def\NN{K}

\begin{document}

\title[A new proof of the atomic decomposition of Hardy spaces]
{A new proof of the atomic decomposition\\ of Hardy spaces}

\author{S. Dekel}
%\author{Shai Dekel}
\address{Hamanofim St. 9, Herzelia, Israel}
\email{Shai.Dekel@ge.com}

\author{G. Kerkyacharian}
%\author{Gerard Kerkyacharian}
\address{
Laboratoire de Probabilit\'{e}s et Mod\`{e}les Al\'{e}atoires, CNRS-UMR 7599,
Universit\'{e} Paris VI et Universit\'{e} Paris VII, rue de Clisson, F-75013 Paris}
\email{kerk@math.univ-paris-diderot.fr}

\author{G. Kyriazis}
%\author{George Kyriazis}
\address{Department of Mathematics and Statistics,
University of Cyprus, 1678 Nicosia, Cyprus}
\email{kyriazis@ucy.ac.cy}

\author{P. Petrushev}
%\author{Pencho Petrushev}
\address{Department of Mathematics\\University of South Carolina\\
Columbia, SC 29208}
\email{pencho@math.sc.edu}



%\thanks{Corresponding author: Pencho Petrushev,
%Email: pencho@math.sc.edu}

\begin{abstract}
A new proof is given of the atomic decomposition of Hardy spaces $H^p$, $0<p\le 1$,
in the classical setting on $\R^n$.
The new method can be used to establish atomic decomposition of maximal Hardy spaces
in general and nonclassical settings.
\end{abstract}

%\date{June 21, 2014}
%\date{July 14, 2014}
%\date{July 15, 2014}

\maketitle

\section{Introduction}\label{sec:introduction}
\setcounter{equation}{0}

The study of the real-variable Hardy spaces $H^p$, $0<p\le 1$,
on $\R^n$ was pioneered by Stein and Weiss \cite{Stein-Weiss}
and a major step forward in developing this theory was made by Fefferman and Stein in \cite{FS}, see also \cite{Stein}.
Since then there has been a great deal of work done on Hardy spaces.
The atomic decomposition of $H^p$ was first established by Coifman \cite{Coifman} in dimension $n=1$
and by Latter \cite{Latter} in dimensions $n>1$.

The purpose of this article is to give a new proof of the atomic decomposition of the $H^p$ spaces %, $0<p\le 1$,
in the classical setting on $\R^n$.
Our method does not use the Calder\'{o}n-Zygmund decomposition of functions
and an approximation of the identity as the classical argument does, see \cite{Stein}.
The main advantage of the new proof over the classical one is that it
is amenable to utilization in more general and nonclassical settings.
For instance, it is used in \cite{DKKP-Hardy} for establishing the equivalence of maximal and atomic Hardy spaces
in the general setting of a metric measure space with the doubling property
and in the presence of a non-negative self-adjoint operator whose heat kernel has Gaussian localization
and the Markov property.

\smallskip

\noindent
{\bf Notation.}
For a set $E\subset \R^n$ we will denote
$E+B(0, \delta):= \cup_{x\in E} B(x, \delta)$, %which is the $\delta$-neighborhood of $E$.
where $B(x, \delta)$ stands for the open ball centered at $x$ of radius $\delta$.
We will also use the notation $cB(x, \delta):= B(x, c\delta)$.
Positive constants will be denoted by $c, c_1, \dots$ and they may vary
at every occurrence; $a\sim b$ will stand for $c_1\le a/b\le c_2$. %, $c_1, c_2>0$.

\subsection{Maximal operators and \boldmath $H^p$ spaces}\label{subsec:maximal-oper}

We begin by recalling some basic facts about Hardy spaces on $\R^n$.
For a complete account of Hardy spaces we refer the reader to \cite{Stein}.

%%%%%%% Definition

Given $\varphi\in \cS$ with $\cS$ being the Schwartz class on $\R^n$ and $f\in \cS'$ one defines
\begin{equation}\label{def-max-1}
M_\varphi f(x):= \sup_{t>0} |\varphi_t*f(x)|
\;\;\hbox{with}\;\;
\varphi_t(x):= t^{-n}\varphi(t^{-1}x),
\quad\hbox{and}
\end{equation}
\begin{equation}\label{def-max-2}
M^*_{\varphi, a} f(x):= \sup_{t>0}\sup_{y\in \R^n, |x-y|\le at} |\varphi_t*f(y)|, \quad a\ge 1.
\end{equation}
We now recall the grand maximal operator.
%
%%%%%%% Definition
%
Write
$$
\fN_N (\varphi):= \sup_{x\in \R^n} (1+|x|)^N \max_{|\alpha|\le N+1}|\partial^\alpha\varphi(x)|
$$
and denote
$$
\cF_N:= \{\varphi\in \cS: \fN_N(\varphi)\le 1\}.
$$
The grand maximal operator is defined by
\begin{equation}\label{def-grand-max1}
\cM_N f(x):= \sup_{\varphi\in\cF_N} M^*_{\varphi, 1}f(x), \quad f\in \cS'.
\end{equation}

It is easy to see that for any $\varphi\in \cS$ and $a\ge 1$ one has
\begin{equation}\label{est-grand-max}
M^*_{\varphi, a} f(x) \le a^N \cP_N(\varphi)\cM_{N} f(x), \quad f\in \cS'.
\end{equation}

%%%%%%% Definition

\begin{definition}\label{def:Hp}
The space $H^p$, $0<p\le 1$, is defined as the set of all bounded distributions
$f\in \cS'$ such that the Poisson maximal function $\sup_{t>0}|P_t*f(x)|$ belongs to $L^p$;
the quasi-norm on $H^p$ is defined by
\begin{equation}\label{def-Hp}
\|f\|_{H^p}:= \big\|\sup_{t>0}|P_t*f(\cdot)|\big\|_{L^p}.
\end{equation}

\end{definition}

As is well known the following assertion holds, see \cite{FS, Stein}:

%%%%%%% Proposition

\begin{proposition}\label{prop:Hp}
Let $0<p\le 1$, $a\ge 1$, and assume $\varphi\in\cS$ and $\int_{\R^n}\varphi\ne 0$. %and $\fN_N(\varphi)\le 1$.
Then for any $N \ge \lfloor\frac{n}{p}\rfloor+1$
\begin{equation}\label{equiv-Hp-norms}
\|f\|_{H^p} \sim \|M^*_{\varphi, a}f\|_{L^p}
\sim \|\cM_N f\|_{L^p},
\quad \forall f\in H^p.
\end{equation}

\end{proposition}

\subsection{Atomic \boldmath $H^p$ spaces}\label{subsec:atomic-Hp-spaces}

%%%%%%%%%% Definition

%\begin{definition} \label{def-atoms}
A~function $a\in L^\infty(\R^n)$ is called an atom if
there exists a ball $B$  such that
\begin{enumerate}
\item[(i)] $\supp a \subset B$,

\medskip

\item[(ii)] $\|a\|_{L^\infty} \le |B|^{-1/p}$, and

\medskip

\item[(iii)] $\int_{\R^n}x^\alpha a(x)dx=0$ for all $\alpha$ with $|\alpha|\le n(p^{-1}-1)$.
\end{enumerate}

%\end{definition}

%%%%%%%% Defintion

%\begin{definition} \label{def-Hardy-Atomic}
The atomic Hardy space $H^p_A$, $0<p\le 1$, is defined as the set of all distributions $f\in \cS'$
that can be represented in the form
\begin{equation}\label{atomic-Hardy}
f=\sum_{j=1}^\infty \lambda_j a_j,
\quad\hbox{where}\quad
\sum_{j=1}^\infty |\lambda_j|^p <\infty,
\end{equation}
$\{a_j\}$ are atoms, and the convergence is in $\cS'$.
Set
\begin{equation}\label{norm-atomic-Hardy}
\|f\|_{H^p_A}:=\inf_{f=\sum_{j} \lambda_j a_j} \Big(\sum_{j=1}^\infty |\lambda_j|^p\Big)^{1/p},
\quad f\in H^p_A.
\end{equation}
%\end{definition}

\section{Atomic decomposition of $H^p$ spaces}\label{sec:equivalence-Hp-spaces}
\setcounter{equation}{0}

We now come to the main point in this article, that is, to give a new proof of the following
classical result \cite{Coifman, Latter}, see also  \cite{Stein}:

%%%%%%%%%%%%% Theorem

\begin{theorem}\label{thm:hardy}
For any  $0<p\le 1$ the continuous embedding $H^p \subset H^p_A$ is valid, that is,
if $f\in H^p$, then $f\in H^p_A$  and
\begin{equation}\label{est-hardy}
\|f\|_{H^p_A} \le c\|f\|_{H^p},
\end{equation}
where $c>0$ is a constant depending only on $p, n$.
This along with the easy to prove embedding  $H^p_A \subset H^p$
leads to $H^p=H^p_A$ and
$\|f\|_{H^p} \sim \|f\|_{H^p_A}$ for $f\in H^p$.
\end{theorem}

\noindent
{\bf Proof.}
We first derive a simple {\em decomposition} identity which will play a central r\^{o}le in this proof.
For this construction we need the following

%%%%%%%% Lemma

\begin{lemma}\label{lem:construct-phi}
For any $m\ge 1$ there exists a function
$\varphi\in C^\infty_0(\R^n)$ such that
$\supp \varphi \subset B(0, 1)$,
$\hat\varphi(0)=1$, and
$\partial^\alpha\hat\varphi(0)=0$ for $0<|\alpha|\le m$.
Here $\hat\varphi(x) :=\int_{\R^n}\varphi(x)e^{-ix\cdot\xi}dx$.
\end{lemma}

\noindent
{\bf Proof.}
We will construct a function $\varphi$ with the claimed properties in dimension $n=1$.
Then a normalized dilation of $\varphi(x_1)\varphi(x_2)\cdots\varphi(x_n)$
will have the claimed properties on $\R^n$.

For the univariate construction, pick a smooth ``bump" $\oomega$ with the following properties:
$\oomega\in C^\infty_0(\R)$, $\supp \oomega \subset [-1/4, 1/4]$,
$\oomega(x)>0$ for $x\in (-1/4, 1/4)$, and $\oomega$ is even.
Let $\Theta(x):= \oomega(x+1/2) - \oomega(x-1/2)$ for $x\in \R$.
Clearly $\Theta$ is odd.

We may assume that $m \ge 1$ is even, otherwise we work with $m+1$ instead.
Denote $\Delta_h^m:= (T_h-T_{-h})^m$, where
%$T_h$ the shift to the left by $h$ operator, defined by
$T_hf(x):= f(x+h)$.

We define
$\varphi(x):= \frac{1}{x}\Delta_h^m\Theta (x)$, where $h=\frac{1}{8m}$.
Clearly, $\varphi \in C^\infty(\R)$,
$\varphi$ is even, and
$\supp \varphi \subset [-\frac{7}{8}, -\frac{1}{8}]\cup [\frac{1}{8}, \frac{7}{8}]$.
It is readily seen that for $\nu=1, 2, \dots, m$
$$
\hat \varphi ^{(\nu)}(\xi)
= (-i)^\nu \int_\R x^{\nu-1}\Delta_h^m\Theta(x)e^{-i\xi x}dx
$$
and hence
$$
\hat \varphi ^{(\nu)}(0)
= (-i)^\nu \int_\R x^{\nu-1}\Delta_h^m\Theta(x)dx
= (-i)^{\nu+m} \int_\R \Theta(x)\Delta_h^m x^{\nu-1}dx =0.
$$
On the other hand,
\begin{align*}
\hat \varphi(0)= \int_\R \varphi(x)dx
=2\int_0^\infty x^{-1}\Delta_h^m\Theta(x)dx
=2(-1)^m\int_{1/4}^{3/4}\Theta(x) \Delta_h^m x^{-1}dx.
\end{align*}
However, for any sufficiently smooth function $f$ we have
$\Delta_h^mf(x) = (2h)^mf^{(m)}(\xi)$, where $\xi \in (x-mh, x+mh)$.
Hence,
$$
\Delta_h^m x^{-1} = (2h)^m m!(-1)^m \xi^{-m-1}
\quad\hbox{with}\quad \xi \in (x-mh, x+mh)\subset [1/8, 7/8].
$$
Consequently, $\hat \varphi(0)\ne 0$ and then $\hat \varphi(0)^{-1}\varphi(x)$
has the claimed properties.
$\qed$

\smallskip

With the aid of the above lemma,
we pick $\varphi\in C^\infty_0(\R^n)$
with the following properties:
$\supp \varphi \subset B(0, 1)$,
$\hat\varphi(0)=1$, and
$\partial^\alpha\hat\varphi(0)=0$ for $0<|\alpha|\le \KK$,
where $K$ is sufficiently large.
More precisely, we choose $\KK \ge n/p$.

Set $\psi(x):= 2^n\varphi(2x)-\varphi(x)$.
Then $\hat\psi(\xi)=\hat\varphi(\xi/2)-\hat\varphi(\xi)$.
Therefore,
$\partial^\alpha\hat\psi(0)=0$ for $|\alpha|\le \KK$
which implies
$\int_{\R^n}x^\alpha \psi(x)dx=0$ for $|\alpha|\le \KK$.
We also introduce the function $\tilde\psi(x):= 2^n\varphi(2x)+\varphi(x)$.
We will use the notation
$h_k(x):= 2^{kn}h(2^kx)$.

Clearly, for any $f\in\cS'$ we have $f=\lim_{j\to\infty}\varphi_j*\varphi_j*f$ (convergence in $\cS'$),
which leads to the following representation:
For any $j\in\bZ$
\begin{align*}
f&=\varphi_j*\varphi_j*f + \sum_{k=j}^\infty \big[\varphi_{k+1}*\varphi_{k+1}*f-\varphi_k*\varphi_k*f\big]\\
&=\varphi_j*\varphi_j*f + \sum_{k=j}^\infty
\big[\varphi_{k+1}-\varphi_k\big]*
\big[\varphi_{k+1}+\varphi_k\big]*f.
\end{align*}
Thus we arrive at
\begin{equation}\label{main-repr}
f=\varphi_j*\varphi_j*f + \sum_{k=j}^\infty \psi_k*\tilde\psi_k*f,
\quad \forall f\in \cS' \;\;\forall j\in\bZ
\quad\hbox{(convergence in $\cS'$).}
\end{equation}
Observe that $\supp \psi_k \subset B(0, 2^{-k})$ and $\supp \tilde\psi_k \subset B(0, 2^{-k})$.

%%%%%%%%%%%%%%%%%%%%%%%%%%%%%

\smallskip

In what follows we will utilize the grand maximal operator $\cM_N$, defined in (\ref{def-grand-max1})
with $N := \lfloor\frac{n}{p}\rfloor+1$.
The following claim follows  readily from (\ref{est-grand-max}):
If $\phi\in \cS$,  % and assume $\fN_N (\phi)\le c$.
then for any $f\in \cS'$, $k\in\bZ$, and $x\in\R^n$
\begin{equation}\label{argument}
|\phi_k*f(y)| \le c\cM_N f(x)
\quad\hbox{for all}\quad y\in\R^n \;\;{with}\;\; |y-x|\le 2^{-k+1},
\end{equation}
where the constant $c>0$ depends only on $\cP_N(\phi)$ and $N$.
%Here as before $\phi_k(x):=2^{kn}\phi(2^kx)$.

\smallskip

Let $f\in H^p$, $0<p\le 1$, $f\ne 0$. We define
\begin{equation}\label{def-Omega-r}
\Omega_r:= \{x\in \R^n: \cM_N f(x) >2^r\}, \quad r\in \bZ.
\end{equation}
Clearly, $\Omega_r$ is open,
$\Omega_{r+1} \subset \Omega_r$, and $\R^n = \cup_{r\in \bZ} \Omega_r$.
%The latter identity follows by $\cM_N f(x) >0$ $\forall x\in \R^n$ due to $f\ne 0$.
%
It is easy to see that %(cf. (5.14) in \cite{DKKP})
\begin{equation}\label{ineq-Omega}
\sum_{r\in \bZ} 2^{pr}|\Omega_r| \le c \int_{\R^n} \cM_N f(x)^p d\mu(x) \le c\|f\|_{H^p}^p.
\end{equation}

%\noindent
%{\bf Convergence.}
From (\ref{ineq-Omega}) we get $|\Omega_r|\le c2^{-pr}\|f\|_{H^p}^p$ for $r\in\bZ$.
Therefore,
%from the definition of $\Omega_r$ in (\ref{def-Omega-r}) and from above, it follows that
for any $r\in\bZ$ there exists $J>0$ such that
%$\|\varphi_j*f\|_\infty \le 2^r$ and hence
$\|\varphi_j*\varphi_j*f\|_\infty \le c2^r$ for $j<-J$.
Consequently, $\|\varphi_j*\varphi_j*f\|_\infty \to 0$ as $j\to -\infty$, which implies
\begin{equation}\label{convergence}
f=\lim_{K\to\infty} \sum_{k=-\infty}^K \psi_k*\tilde\psi_k*f
\quad \hbox{(convergence in $\cS'$).}
%\quad \hbox{in} \;\;\cS'.
\end{equation}

Assuming that $\Omega_r\ne \emptyset$ we write
$$ %\begin{equation}\label{def-E-rk}
E_{rk}:=\big\{x\in \Omega_r: \dist(x, \Omega_r^c) > 2^{-k+1}\big\}
\setminus\big\{x\in \Omega_{r+1}: \dist(x, \Omega_{r+1}^c) > 2^{-k+1}\big\}.
$$ %\end{equation}
By (\ref{ineq-Omega}) it follows that $|\Omega_r| <\infty$ and hence there exists $\kr\in \bZ$
such that $E_{r\kr}\ne\emptyset$ and $E_{rk}=\emptyset$ for $k<\kr$.
Evidently $s_r \le s_{r+1}$.
We define
\begin{equation}\label{def-F-rr}
F_{r}(x):= \sum_{k \ge \kr} \int_{E_{rk}} \psi_k(x-y)\tilde\psi_k*f(y) dy,
\quad x\in \R^n, \; r\in \bZ,
\end{equation}
and more generally
\begin{equation}\label{def-F-rk}
F_{r, \kl, \kll}(x):= \sum_{k = \kl}^{\kll} \int_{E_{rk}} \psi_k(x-y)\tilde\psi_k*f(y) dy,
\quad \kr\le \kl \le \kll\le \infty.
\end{equation}
It will be shown in Lemma~\ref{lem:F-rr} below that the functions $F_r$ and $F_{r, \kl, \kll}$
are well defined and $F_r, F_{r, \kl, \kll}\in L^\infty$.

Note that $\supp \psi_k \subset B(0, 2^{-k})$ and hence
\begin{equation}\label{supp-block}
\supp \Big(\int_{E_{rk}} \psi_k(x-y)\tilde\psi_k*f(y) dy\Big) \subset E_{rk}+B(0,2^{-k}).
\end{equation}
On the other hand, clearly
$2B(y,2^{-k}) \cap \big(\Omega_r\setminus\Omega_{r+1}\big) \ne \emptyset$
for each $y\in E_{rk}$,
and $\fN_N(\tilde\psi) \le c$.
Therefore, see (\ref{argument}),
$|\tilde\psi_k*f(y)|\le c2^r$ for $y\in E_{rk}$, which implies
\begin{equation}\label{norm-block}
\Big\|\int_{E} \psi_k(\cdot-y)\tilde\psi_k*f(y) dy\Big\|_\infty \le c2^r,\quad \forall E\subset E_{rk}.
\end{equation}
Similarly,
\begin{equation}\label{norm-block-2}
\Big\|\int_{E} \varphi_k(\cdot-y)\tilde\varphi_k*f(y) dy\Big\|_\infty \le c2^r, \quad \forall E\subset E_{rk}.
\end{equation}

We collect all we need about the functions $F_r$ and $F_{r, \kl, \kll}$ in the following

%%%%%%%% Lemma

\begin{lemma}\label{lem:F-rr}
$(a)$
We have
\begin{equation}\label{Erk-cover}
E_{rk}\cap E_{r'k}=\emptyset
\quad\hbox{if $r\ne r'\;\;$ and}
\quad\R^n=\cup_{r\in \bZ} E_{rk}, \quad \forall k\in\bZ.
\end{equation}

$(b)$
There exists a constant $c>0$ such that for any $r\in\bZ$ and $s_r\le \kl \le \kll\le \infty$
\begin{equation}\label{est-F-rr}
\|F_r\|_\infty \le c2^r, \quad \|F_{r, \kl, \kll}\|_\infty \le c2^r.
\end{equation}

$(c)$
The series in $(\ref{def-F-rr})$ and $(\ref{def-F-rk})$ $($if $\kll=\infty$$)$ converge point-wise
and in distributional sense.

$(d)$
Moreover,
\begin{equation}\label{F-rr-zero}
F_r(x)=0, %\quad\hbox{for}\;\;
\quad\forall x\in \R^n\setminus\Omega_r, \;\; \forall r\in\bZ.
\end{equation}
\end{lemma}

\noindent
{\bf Proof.}
Identities (\ref{Erk-cover}) are obvious
and (\ref{F-rr-zero}) follows readily from (\ref{supp-block}).

We next prove the left-hand side inequality in (\ref{est-F-rr});
the proof of the right-hand side inequality is similar and will be omitted.
Consider the case when $\Omega_{r+1}\ne \emptyset$
(the case when $\Omega_{r+1}= \emptyset$ is easier).
Write
$$ %\begin{equation}\label{def-U-V}
U_k=\big\{x\in \Omega_r: \dist(x, \Omega_r^c) > 2^{-k+1}\big\},
\quad
V_k=\big\{x\in \Omega_{r+1}: \dist(x, \Omega_{r+1}^c) > 2^{-k+1}\big\}.
$$ %\end{equation}
Observe that $E_{rk} = U_k\setminus V_k$.

From (\ref{supp-block}) it follows that $|F_r(x)|=0$ for
$x \in \R^n\setminus \cup_{k\ge \kr}(E_{rk}+B(0, 2^{-k}))$.
We next estimate $|F_r(x)|$ for $x \in \cup_{k\ge \kr}(E_{rk}+B(0, 2^{-k}))$.
Two cases present themselves here.

{\bf Case 1:}
$x\in \big[\cup_{k\ge \kr} (E_{rk}+B(0,2^{-k}))\big]\cap\Omega_{r+1}$.
Then there exist $\nu, \ell\in\bZ$ such that
\begin{equation}\label{location-x}
x\in (U_{\ell+1}\setminus U_\ell)\cap (V_{\nu+1}\setminus V_\nu).
\end{equation}
Due to $\Omega_{r+1}\subset\Omega_r$ we have $V_k\subset U_k$, implying
$(U_{\ell+1}\setminus U_\ell)\cap (V_{\nu+1}\setminus V_\nu)=\emptyset$ if $\nu<\ell$.
We consider two subcases depending on whether $\nu\ge \ell+3$ or $\ell \le \nu \le \ell+2$.

(a) Let $\nu\ge \ell+3$.
We claim that (\ref{location-x}) yields
\begin{equation}\label{intersect-1}
B(x, 2^{-k}) \cap E_{rk} =\emptyset \quad\hbox{for}\quad k\ge \nu+2 \;\;\hbox{or}\;\; k\le \ell-1.
\end{equation}
Indeed, if $k\ge \nu+2$, then $E_{rk} \subset \Omega_r\setminus V_{\nu+2}$, which implies (\ref{intersect-1}),
while if $k\le \ell-1$, then $E_{rk} \subset U_{\ell-1}$, again implying (\ref{intersect-1}).

We also claim that
\begin{equation}\label{claim-1}
B(x, 2^{-k}) \subset E_{rk}\quad\hbox{for}\quad \ell+2\le k\le \nu-1.
\end{equation}
Indeed, clearly
\begin{align*}
(U_{\ell+1}\setminus U_\ell)\cap (V_{\nu+1}\setminus V_\nu)
\subset (U_{k-1}\setminus U_\ell)\cap (V_{\nu+1}\setminus V_{\nu+1})
\subset U_{k-1}\setminus V_{k+1},
\end{align*}
which implies (\ref{claim-1}).

From (\ref{supp-block}) and (\ref{intersect-1})- (\ref{claim-1}) it follows that
\begin{align*}
F_{r}(x)&= \sum_{k =\ell}^{\nu+1} \int_{E_{rk}} \psi_k(x-y)\tilde\psi_k*f(y) dy
= \sum_{k =\ell}^{\ell+1} \int_{E_{rk}} \psi_k(x-y)\tilde\psi_k*f(y) dy\\
&+ \sum_{k =\ell+2}^{\nu-2} \int_{\R^n} \psi_k(x-y)\tilde\psi_k*f(y) dy
+ \sum_{k =\nu-1}^{\nu+1} \int_{E_{rk}} \psi_k(x-y)\tilde\psi_k*f(y) dy.
\end{align*}
However,
\begin{align*}
&\sum_{k =\ell+2}^{\nu-2} \int_{\R^n} \psi_k(x-y)\tilde\psi_k*f(y) dy
%= \sum_{k =\ell+2}^{\nu-1}\psi_k*\tilde\psi_k*f(x)\\
%&\qquad\qquad
= \sum_{k =\ell+2}^{\nu-2} \big[\varphi_{k+1}*\varphi_{k+1}*f(x)-\varphi_k*\varphi_k*f(x)\big]\\
&\qquad\qquad
= \varphi_{\nu-1}*\varphi_{\nu-1}*f(x)-\varphi_{\ell+2}*\varphi_{\ell+2}*f(x)\\
&\qquad\qquad
=\int_{E_{r, \nu-1}} \varphi_{\nu-1}(x-y)\varphi_{\nu-1}*f(y) dy
- \int_{E_{r, \ell+2}} \varphi_{\ell+2}(x-y)\varphi_{\ell+2}*f(y) dy.
\end{align*}
Combining the above with (\ref{norm-block}) and (\ref{norm-block-2}) we obtain
$|F_r(x)| \le c2^r$.

(b) Let $\ell \le \nu \le \ell+2$.
Just as above we have
\begin{align*}
F_{r}(x)&= \sum_{k =\ell}^{\nu+1} \int_{E_{rk}} \psi_k(x-y)\tilde\psi_k*f(y) dy
= \sum_{k =\ell}^{\ell+3} \int_{E_{rk}} \psi_k(x-y)\tilde\psi_k*f(y) dy
\end{align*}
We use (\ref{norm-block}) to estimate each of these four integrals and again obtain
$|F_r(x)| \le c2^r$.

\smallskip

%%%%%%%%%%%%%%%%%%%%%%%%%

{\bf Case 2:} $x\in \Omega_r\setminus\Omega_{r+1}$.
Then there exists $\ell \ge \kr$ such that
$$
x\in (U_{\ell+1}\setminus U_\ell)\cap(\Omega_r\setminus\Omega_{r+1}).
$$
Just as in the proof of (\ref{intersect-1}) we have
$B(x, 2^{-k}) \cap E_{rk} =\emptyset$ for $k\le \ell-1$,
and as in the proof of (\ref{claim-1}) we have
$$ %\begin{align*}
(U_{\ell+1}\setminus U_\ell)\cap (\Omega_r\setminus \Omega_{r+1})
\subset U_{k-1}\setminus V_{k+1},
$$ %\end{align*}
which implies
%$ %\begin{equation}\label{claim-3}
%B(x, 2^{-k}) \subset E_{rk}\quad\hbox{for}\quad k\ge \ell+2.
%$ %\end{equation}
$B(x, 2^{-k}) \subset E_{rk}$ for $k\ge \ell+2$.
We use these and (\ref{supp-block}) to obtain
\begin{align*}
F_{r}(x)
&= \sum_{k =\ell}^\infty \int_{E_{rk}} \psi_k(x-y)\tilde\psi_k*f(y) dy\\
&= \sum_{k =\ell}^{\ell+1} \int_{E_{rk}} \psi_k(x-y)\tilde\psi_k*f(y) dy
+ \sum_{k =\ell+2}^\infty \int_{\R^n} \psi_k(x-y)\tilde\psi_k*f(y) dy.
\end{align*}
For the last sum we have
\begin{align*}
&\sum_{k =\ell+2}^\infty \int_{\R^n} \psi_k(x-y)\tilde\psi_k*f(y) dy
= \lim_{\nu\to \infty}\sum_{k =\ell+2}^{\nu}\psi_k*\tilde\psi_k*f(x)\\
&= \lim_{\nu\to \infty}\big(\varphi_{\nu+1}*\varphi_{\nu+1}*f(x)-\varphi_{\ell+2}*\varphi_{\ell+2}*f(x))\\
&= \lim_{\nu\to \infty}\Big(\int_{E_{r, \nu+1}} \varphi_{\nu+1}(x-y)\varphi_{\nu+1}*f(y) dy
- \int_{E_{r, \ell+2}} \varphi_{\ell+2}(x-y)\varphi_{\ell+2}*f(y) dy\Big).
\end{align*}
From the above and (\ref{norm-block})-(\ref{norm-block-2}) we obtain
$|F_r(x)| \le c2^r$.

The point-wise convergence of the series in (\ref{def-F-rr}) follows from above
and we similarly establish the point-wise convergence in (\ref{def-F-rk}).

The convergence in distributional sense in (\ref{def-F-rr}) relies on the following assertion:
For every $\phi\in\cS$
\begin{equation}\label{conv-distr}
\sum_{k \ge \kr} |\langle g_{rk}, \phi\rangle| <\infty,
\quad\hbox{where}\quad g_{rk}(x):=\int_{E_{rk}} \psi_k(x-y)\tilde\psi_k*f(y) dy.
\end{equation}
Here $\langle g_{rk}, \phi\rangle:= \int_{\R^n} g_{rk} \overline{\phi} dx$.
To prove the above we will employ this estimate:
\begin{equation}\label{est-psi-f}
\|\tilde\psi_k f\|_\infty \le c2^{kn/p}\|f\|_{H^p}, \quad k\in \bZ.
\end{equation}
Indeed, using (\ref{est-grand-max}) we get
\begin{align*}
|\tilde\psi_k f(x)|^p
&\le \inf_{y: |x-y|\le 2^{-k}} \sup_{z: |y-z|\le 2^{-k}}|\tilde\psi_k f(z)|^p
\le \inf_{y: |x-y|\le 2^{-k}} c\cM_N(f)(y)^p\\
& \le c|B(x, 2^{-k})|^{-1}\int_{B(x, 2^{-k})}\cM_N(f)(y)^p d\mu(y)
\le c2^{kn}\|f\|_{H^p}^p,
\end{align*}
and (\ref{est-psi-f}) follows.

We will also need the following estimate:
For any $\sigma>n$ there exists a constant $c_\sigma>0$ such that
\begin{equation}\label{inner-prod}
\Big|\int_{\R^n} \psi_k(x-y)\phi(x) dx\Big| \le c_\sigma 2^{-k(K+1)}(1+|y|)^{-\sigma},
\quad y\in \R^n, \; k\ge 0.
\end{equation}
This is a standard estimate for inner products taking into account that
$\phi\in\cS$ and $\psi\in C^\infty$, $\supp \psi \subset B(0, 1)$, and
$\int_{\R^n}x^\alpha\psi(x)dx=0$ for $|\alpha|\le K$.

We now estimate $|\langle g_{rk}, \phi\rangle|$.
From (\ref{est-psi-f}) and the fact that $\psi\in C^\infty_0(\R)$ and $\phi\in\cS$ it readily follows that
$$
\int_{E_{rk}}\int_{\R^n} |\psi_k(x-y)||\phi(x)||\tilde\psi_kf(y)| dy dx <\infty,
\quad k \ge s_r.
$$
Therefore, we can use Fubini's theorem, (\ref{est-psi-f}), and (\ref{inner-prod}) to obtain for $k\ge 0$
\begin{align}\label{gk-phi}
|\langle g_{rk}, \phi\rangle|
& \le \int_{E_{rk}} \Big|\int_{\R^n} \psi_k(x-y)\phi(x) dx\Big||\tilde\psi_kf(y)| dy\notag\\
%& = \int_{E_{rk}} \Big|\int_{\R^n} \psi_k(y, x)\overline{\phi(x)} dx\Big||\tilde\psi_kf(y)| dy\\
& \le c2^{-k(\NN+1-n/p)}\|f\|_{H^p} \int_{E_{rk}}(1+|y|)^{-\sigma}dy
\le c2^{-k(\NN+1-n/p)}\|f\|_{H^p},
\end{align}
which implies (\ref{conv-distr}) because $\NN \ge n/p$.

Denote $G_\ell := \sum_{k=s_r}^\ell g_{rk}$.
From the above proof of (b) and (\ref{est-F-rr}) we infer that $G_\ell(x) \to F_r(x)$ as $\ell \to \infty$
for $x\in \R^n$ and $\|G_\ell\|_\infty \le c2^r< \infty$ for $\ell \ge s_r$.
On the other hand, from (\ref{conv-distr}) it follows that the series $\sum_{k\ge s_r} g_{rk}$
converges in distributional sense.
By applying the dominated convergence theorem one easily concludes that $F_r=\sum_{k\ge s_r} g_{rk}$
with the convergence in distributional sense.
$\qed$

\smallskip

We set $F_r:= 0$ in the case when $\Omega_r=\emptyset$, $r\in\bZ$.

Note that by (\ref{Erk-cover}) it follows that
\begin{equation}\label{rep-psi-psi}
\psi_k*\psi_k*f(x)= \int_{\R^n}\psi_k(x-y)\psi_k*f(y)dy
=\sum_{r\in\bZ}\int_{E_{rk}} \psi_k(x-y)\tilde\psi_k*f(y) dy
\end{equation}
and using $(\ref{convergence})$ and the definition of $F_r$ in $(\ref{def-F-rr})$ we arrive at
\begin{equation}\label{represent}
f = \sum_{r\in\bZ} F_r \;\; \hbox{in}\;\; \cS',
\;\;\hbox{i.e.}\quad
\langle f, \phi\rangle = \sum_{r\in\bZ} \langle F_r, \phi\rangle,
\quad\forall \phi\in\cS,
\end{equation}
where the last series converges absolutely.
Above $\langle f, \phi\rangle$ denotes the action of $f$ on~$\overline{\phi}$.
We next provide the needed justification of identity (\ref{represent}).

From (\ref{convergence}), (\ref{def-F-rr}), (\ref{rep-psi-psi}),
and the notation from (\ref{conv-distr}) we obtain for $\phi\in\cS$
\begin{align*}
\langle f, \phi\rangle
= \sum_{k}\langle\psi_k\tilde\psi_kf, \phi \rangle
= \sum_{k}\sum_{r}\langle g_{rk}, \phi\rangle
= \sum_{r}\sum_{k}\langle g_{rk}, \phi\rangle
= \sum_{r}\langle F_r, \phi\rangle.
\end{align*}
Clearly, to justify the above identities it suffices to show that
$
\sum_{k}\sum_{r} |\langle g_{rk}, \phi\rangle| <\infty.
$
We split this sum into two:
$
\sum_{k}\sum_{r}\cdots
= \sum_{k\ge 0}\sum_{r} \cdots + \sum_{k<0}\sum_{r} \cdots
=: \Sigma _1+ \Sigma_2.
$
To estimate $\Sigma_1$ we use (\ref{gk-phi}) and obtain
\begin{align*}
\Sigma_1
& \le c\|f\|_{H^p} \sum_{k\ge 0}2^{-k(\NN+1-n/p)} \sum_r\int_{E_{rk}}(1+|y|)^{-\sigma}dy\\
& \le c\|f\|_{H^p} \sum_{k\ge 0}2^{-k(\NN+1-n/p)} \int_{\R^n}(1+|y|)^{-\sigma}dy
\le c\|f\|_{H^p}.
\end{align*}
Here we also used that $K\ge n/p$ and $\sigma >n$.

We estimate $\Sigma_2$ in a similar manner, using
the fact that $\int_{\R^n}|\psi_k(y)|dy \le c<\infty$ and (\ref{est-psi-f}).
We get
\begin{align*}
\Sigma_2
&\le c\|f\|_{H^p}  \sum_{k<0}2^{kn/p}\sum_{r}
\int_{E_{rk}}\int_{\R^n} |\psi_k(x-y)|dy|\phi(x)| dx\\
&\le c\|f\|_{H^p}  \sum_{k<0}2^{kn/p}\int_{\R^n} (1+|x|)^{-n-1} dx
\le c\|f\|_{H^p}.
\end{align*}
The above estimates of $\Sigma_1$ and $\Sigma_2$ imply
$
\sum_{k}\sum_{r} |\langle g_{rk}, \phi\rangle| <\infty,
$
which completes the justification of (\ref{represent}).

\smallskip

Observe that due to $\int_{\R^n}x^\alpha \psi(x)dx=0$ for $|\alpha|\le \NN$ we have
\begin{equation}\label{moments}
\int_{\R^n}x^\alpha F_r(x)dx=0 \quad
\hbox{for $|\alpha|\le \NN$, $r\in\bZ$.}
\end{equation}

%%%%%%%%%%%%%%%%%

We next decompose each function $F_r$ into atoms.
To this end we need a Whitney type cover for $\Omega_r$, given in the following

%%%%%%%% Lemma

\begin{lemma}\label{lem:Whitney}
Suppose $\Omega$ is an open proper subset of $\R^n$ and
let $\rho(x):= \dist (x, \Omega^c)$.
Then there exists a constant $K >0$, depending only on $n$,
and a sequence of points $\{\xi_j\}_{j\in\bN}$ in $\Omega$
with the following properties,
where $\rho_j:= \dist (\xi_j, \Omega^c)$:

\medskip

$(a)$ $\Omega = \cup_{j\in \bN} B(\xi_j, \rho_j/2)$.

\medskip

$(b)$ $\{B(\xi_j, \rho_j/5)\}$ are disjoint.

\medskip

$(c)$ If $B\big(\xi_j, \frac{3\rho_j}{4}\big)\cap B\big(\xi_\nu, \frac{3\rho_\nu}{4}\big)\ne \emptyset$,
then  $7^{-1}\rho_\nu\le \rho_j \le 7\rho_\nu$.

\medskip

$(d)$ For every $j\in\bN$ there are at most $K$ balls $B\big(\xi_\nu, \frac{3\rho_\nu}{4}\big)$
intersecting $B\big(\xi_j, \frac{3\rho_j}{4}\big)$.

\end{lemma}

Variants of this simple lemma are well known and frequently used.
To prove it one simply selects $\{B(\xi_j, \rho(\xi_j)/5)\}_{j\in\bN}$
to be a maximal disjoint subcollection of $\{B(x, \rho(x)/5)\}_{x\in\Omega}$
and then properties (a)-(d) follow readily, see \cite{Stein}, pp. 15-16.
%For completeness we give its proof in the appendix.

%%%%%%%%%%%%%%%%%%%%%%%%%%

\smallskip

We apply Lemma~\ref{lem:Whitney} to each set $\Omega_r\ne \emptyset$, $r\in\bZ$.
Fix $r\in \bZ$ and assume $\Omega_r\ne \emptyset$.
Denote by $B_j:= B(\xi_j, \rho_j/2)$, $j=1, 2, \dots$,
the balls given by Lemma~\ref{lem:Whitney}, applied to $\Omega_r$,
with the additional assumption that these balls are ordered so that
$\rho_1 \ge \rho_2 \ge \cdots$.
We will adhere to the notation from Lemma~\ref{lem:Whitney}.
We will also use the more compact notation
$\cB_r:=\{B_j\}_{j\in\bN}$ for the set of balls covering $\Omega_r$.

For each ball $B\in \cB_r$ and $k\ge \kr$ we define
\begin{equation}\label{def-E-B}
E_{rk}^B:=E_{rk} \cap\big(B+2B(0, 2^{-k})\big)
\quad\hbox{if}\quad B\cap E_{rk}\ne \emptyset
\end{equation}
and set
$E_{rk}^B:=\emptyset$ if $B\cap E_{rk}= \emptyset$.

We also define, for $\ell=1, 2, \dots$,
\begin{equation}\label{def-R-B}
R_{rk}^{B_\ell}:= E_{rk}^{B_\ell}\setminus \cup_{\nu>\ell}E_{rk}^{B_\nu}
\quad \hbox{and} %\ell \ge 1,
\end{equation}
\begin{equation}\label{def-F-B}
F_{B_\ell}(x):= \sum_{k \ge \kr} \int_{R_{rk}^{B_\ell}}
\psi_k(x-y)\tilde\psi_k*f(y) dy.
\end{equation}

%%%%%%%% Lemma

\begin{lemma}\label{lem:FrkB}
For every $\ell\ge 1$ the function $F_{B_\ell}$ is well defined, more precisely,
the series in $(\ref{def-F-B})$ converges point-wise and in distributional sense. % $($in $\cS'$$)$.
Furthermore,
\begin{equation}\label{supp-FB}
\supp F_{B_\ell} \subset 7B_\ell,
\end{equation}
\begin{equation}\label{vanish-mom}
\int_{\R^n}x^\alpha F_{B_\ell}(x)dx=0 \quad\hbox{for all $\alpha$ with $|\alpha|\le n(p^{-1}-1)$,}
\end{equation}
and
\begin{equation}\label{est-FB}
\|F_{B_\ell}\|_\infty \le \cs 2^r,
\end{equation}
where the constant $\cs$ is independent of $r, \ell$.

In addition, for any $k\ge \kr$
\begin{equation}\label{FrkB1}
E_{rk}=\cup_{\ell\ge 1} R_{rk}^{B_\ell}
\quad\hbox{and}\quad R_{rk}^{B_\ell}\cap R_{rk}^{B_m} = \emptyset, \quad \ell\ne m.
\end{equation}
Hence
\begin{equation}\label{decomp-Fr}
F_r= \sum_{B\in\cB_r} F_{B} \quad\hbox{$($convergence in $\cS'$$)$.}
\end{equation}
\end{lemma}

\noindent
{\bf Proof.}
Fix $\ell \ge 1$. Observe that using Lemma~\ref{lem:Whitney} we have
$B_\ell \subset \Omega_r^c + B(0, 2\rho_\ell)$
and hence $E_{rk}^{B_\ell}:=\emptyset$ if $2^{-k+1} \ge 2\rho_\ell$.
Define $\kkl: = \min \{k: 2^{-k} < \rho_\ell\}$.
Hence $\rho_\ell/2 \le 2^{-\kkl} <\rho_\ell$.
Consequently,
\begin{equation}\label{def-F-B-2}
F_{B_\ell}(x):= \sum_{k \ge \kkl} \int_{R_{rk}^{B_\ell}}
\psi_k(x-y)\tilde\psi_k*f(y) dy.
\end{equation}
It follows that
$\supp F_{B_\ell} \subset B\big(\xi_\ell, (7/2)\rho_\ell\big)=7B_\ell$,
which confirms (\ref{supp-FB}).

To prove (\ref{est-FB}) we will use the following

%%%%%%%% Lemma

\begin{lemma}\label{lem:FS}
For an arbitrary set $S\subset \R^n$ let
%$
%S_k:= S+2B(0,2^{-k}),
%$
%which is the same as
$
S_k:= \{x\in \R^n: \dist(x, S) < 2^{-k+1}\}
$
and set
\begin{equation}\label{def-FS}
F_S(x):= \sum_{k \ge \kl} \int_{E_{rk}\cap S_k} \psi_k(x-y)\tilde\psi_k*f(y) dy
\end{equation}
for some $\kl \ge \kr$.
Then
$
\|F_S\|_\infty \le c2^r,
$
where $c>0$ is a constant independent of $S$ and $\kl$.
Moreover, the above series converges in $\cS'$.
\end{lemma}

\noindent
{\bf Proof.}
From (\ref{supp-block}) it follows that $F_S(x)=0$ if $\dist (x, S) \ge 3\times 2^{-\kl}$

Let $x\in S$. Evidently,
$B(x, 2^{-k})\subset S_k$ for every $k$ and hence
\begin{align*}
F_S(x)
&=\sum_{k \ge \kl} \int_{E_{rk}\cap B(x, 2^{-k})} \psi_k(x-y)\tilde\psi_k*f(y) dy\\
&=\sum_{k \ge \kl} \int_{E_{rk}} \psi_k(x-y)\tilde\psi_k*f(y) dy= F_{r,\kl}(x).
\end{align*}
On account of Lemma~\ref{lem:F-rr} (b) we obtain
$|F_S(x)|= |F_{r,\kl}(x)|\le c2^r$.

Consider the case when $x\in S_\ell\setminus S_{\ell+1}$ for some $\ell \ge \kl$.
Then $B(x, 2^{-k}) \subset S_k$ if $\kl\le k\le \ell-1$
and $B(x, 2^{-k}) \cap S_k=\emptyset$ if $k\ge \ell+2$.
Therefore,
\begin{align*}
F_S(x) &=\sum_{k = \kl}^{\ell-1} \int_{E_{rk}} \psi_k(x-y)\tilde\psi_k*f(y) dy
+ \sum_{k = \ell}^{\ell+1} \int_{E_{rk}\cap S_k} \psi_k(x-y)\tilde\psi_k*f(y) dy\\
&= F_{r, \kl, \ell-1}(x)+ \sum_{k = \ell}^{\ell+1} \int_{E_{rk}\cap S_k} \psi_k(x-y)\tilde\psi_k*f(y) dy,
\end{align*}
where we used the notation from (\ref{def-F-rk}).
By Lemma~\ref{lem:F-rr} (b) and (\ref{norm-block})
it follows that $|F_S(x)|\le c2^r$.

We finally consider the case when $2^{-\kl+1} \le \dist (x, S) < 3\times 2^{-\kl}$.
Then we have
$
F_S(x) = \int_{E_{r\kl}\cap S_\kl} \psi_\kl(x-y)\tilde\psi_\kl*f(y) dy
$
and the estimate $|F_S(x)|\le c2^r$ is immediate from (\ref{norm-block}).

The convergence in $\cS'$ in (\ref{def-FS}) is established as in the proof of Lemma~\ref{lem:F-rr}.
$\qed$

\medskip
Fix $\ell \ge 1$ and let $\{B_j: j\in \cJ\}$ be the set of all balls $B_j= B(\xi_j, \rho_j/2)$
such that $j>\ell$ and
$$
B\Big(\xi_j, \frac{3\rho_j}{4}\Big) \cap B\Big(\xi_\ell, \frac{3\rho_\ell}{4}\Big) \ne \emptyset.
$$
By Lemma~\ref{lem:Whitney} it follows that $\# \cJ \le K$
and $7^{-1} \rho_\ell \le \rho_j \le 7\rho_\ell$ for $j\in\cJ$.
Define
\begin{equation}\label{def-k1}
k_1:= \min\Big\{k: 2^{-k+1} < 4^{-1}\min\big\{\rho_j: j\in\cJ \cup\{\ell\}\big\} \Big\}.
\end{equation}
From this definition and  $2^{-\kkl} <\rho_\ell$ we infer
\begin{equation}\label{k-kkl}
2^{-k_1+1} \ge 8^{-1} \min\big\{\rho_j: j\in\cJ \cup\{\ell\}\big\}
> 8^{-2}\rho_\ell > 8^{-2}2^{-\kkl}
\; \Longrightarrow \; k_1 \le \kkl+7.
\end{equation}
%implying $k_1 < \kkl+7$.
Clearly, from (\ref{def-k1})
\begin{equation}\label{BB-1}
B_j + 2B(0,2^{-k}) \subset B\big(\xi_j, 3\rho_j/4\big),
%\quad\hbox{for}
\quad \forall k\ge k_1, \;\; \forall j\in \cJ \cup\{\ell\}.
\end{equation}
Denote
$S:= \cup_{j\in\cJ} B_j$ and
$\tilde{S}:= \cup_{j\in\cJ} B_j\cup B_\ell = S\cup B_\ell$.
As in Lemma~\ref{lem:FS} we set
$$
S_k:= S+2B(0, 2^{-k})
\quad\hbox{and}\quad
\tilde{S}_k:= \tilde{S}+2B(0, 2^{-k}).
$$
It readily follows from the definition of $k_1$ in (\ref{def-k1}) that
\begin{equation}\label{R-EE}
R_{rk}^{B_\ell} := E_{rk}^{B_\ell}\setminus \cup_{\nu>\ell} E_{rk}^{B_\nu}
= \big(E_{rk}\cap \tilde{S}_k\big) \setminus \big(E_{rk}\cap S_k\big)
\quad\hbox{for}\quad k\ge k_1.
\end{equation}
Denote
\begin{align*}
F_S(x)&:= \sum_{k \ge k_1} \int_{E_{rk}\cap S_k} \psi_k(x-y)\tilde\psi_k*f(y) dy,\quad\hbox{and}\\
F_{\tilde{S}}(x)&:= \sum_{k \ge k_1} \int_{E_{rk}\cap \tilde{S}_k} \psi_k(x-y)\tilde\psi_k*f(y) dy.
\end{align*}
From (\ref{R-EE}) and the fact that $S \subset \tilde{S}$ it follows that
\begin{align*}
F_{B_\ell}(x) = F_{\tilde{S}}(x) - F_S(x)
+ \sum_{\kkl\le k < k_1} \int_{R_{rk}^{B_\ell}} \psi_k(x-y)\tilde\psi_k*f(y) dy.
\end{align*}
By Lemma~\ref{lem:FS} we get
$\|F_S\|_\infty \le c2^r$ and $\|F_{\tilde{S}}\|_\infty \le c2^r$.
On the other hand from (\ref{k-kkl}) we have $k_1-\kkl \le 7$.
We estimate each of the (at most $7$) integrals above using (\ref{norm-block}) to conclude that
$\|F_{B_\ell}\|_\infty \le c2^r$.

We deal with the convergence in (\ref{def-F-B}) and (\ref{decomp-Fr})
as in the proof of Lemma~\ref{lem:F-rr}.

Clearly, (\ref{vanish-mom}) follows from the fact that
$\int_{\R^n}x^\alpha \psi(x)dx=0$ for all $\alpha$ with $|\alpha|\le \NN$.

Finally, from Lemma~\ref{lem:Whitney} we have $\Omega_r\subset \cup_{j\in\bN} B_\ell$
and then (\ref{FrkB1}) is immediate from (\ref{def-E-B}) and (\ref{def-R-B}).
$\qed$

%\smallskip
%%%%%%%%%%%%%%%%%%%%%

We are now prepared to complete the proof of Theorem~\ref{thm:hardy}.
For every ball $B\in \cB_r$, $r\in\bZ$, provided $\Omega_r\ne \emptyset$, we define
$B^\star:= 7B$,
$$
a_B(x):= \cs^{-1}|B^\star|^{-1/p}2^{-r}F_{B}(x)
\quad\hbox{and}\quad
\lambda_B:= \cs|B^\star|^{1/p}2^r,
$$
%and
%$\lambda_B:= \cs|B^\star|^{1/p}2^r$,
where $\cs>0$ is the constant
from (\ref{est-FB}).
By (\ref{supp-FB}) $\supp a_B\subset B^\star$
and by (\ref{est-FB})
$$
\|a_B\|_\infty \le \cs^{-1}|B^\star|^{-1/p}2^{-r}\|F_{B}\|_\infty \le |B^\star|^{-1/p}.
$$
Furthermore, from (\ref{vanish-mom}) it follows that
$\int_{\R^n}x^\alpha a_B(x)dx=0$ if $|\alpha|\le n(p^{-1}-1)$.
Therefore, each $a_B$ is an atom for $H^p$.

We set $\cB_r:=\emptyset$ if $\Omega_r= \emptyset$.
Now, using the above, (\ref{represent}), and Lemma~\ref{lem:FrkB} we get

$$
f=\sum_{r\in\bZ} F_r
=\sum_{r\in\bZ}\sum_{B\in \cB_r}F_{B}
=\sum_{r\in\bZ}\sum_{B\in \cB_r}\lambda_B a_B,
$$
where the convergence is in $\cS'$, and
$$
\sum_{r\in\bZ}\sum_{B\in \cB_r}|\lambda_B|^p
\le c\sum_{r\in\bZ} 2^{pr}\sum_{B\in \cB_r}|B|
= c\sum_{r\in\bZ} 2^{pr}|\Omega_r|
\le c \|f\|_{H^p}^p,
$$
which is the claimed atomic decomposition of $f\in H^p$.
Above we used that $|B^\star| = |7B| = 7^n|B|$.
$\qed$

%%%%%%%%% remark

\begin{remark}\label{rem:q-atoms}
The proof of Theorem~\ref{thm:hardy} can be considerably simplified and shortened if
one seeks to establish atomic decomposition of the $H^p$ spaces in terms of
$q$-atoms with $p<q<\infty$ rather than $\infty$-atoms as in Theorem~\ref{thm:hardy},
i.e. atoms satisfying
$\|a\|_{L^q} \le |B|^{1/q-1/p}$ with $q<\infty$
rather than
$\|a\|_{L^\infty} \le |B|^{-1/p}$.
We will not elaborate on this here.
\end{remark}


\begin{thebibliography}{999999}

\bibitem{Coifman}
    R. Coifman, A real variable characterization of $H^p$, Studia Math. 51 (1974), 269–-274.

%\bibitem{DKKP}
%    S. Dekel, G, Kyriazis, G. Kerkyacharian, P. Petrushev,
%    Compactly supported frames for spaces of distributions in the framework of Dirichlet spaces, preprint

\bibitem{DKKP-Hardy}
    S. Dekel, G, Kyriazis, G. Kerkyacharian, P. Petrushev,
    Hardy spaces associated with non-negative self-adjoint operators, preprint.
    %Hardy spaces in the general framework of Dirichlet spaces, preprint.

\bibitem{FS}
    C. Fefferman, E. Stein, $H^p$ spaces of several variables, Acta Math. 129 (1972), 137-–193.

\bibitem{Latter}
    R. Latter, A characterization of $H^p(\R^n)$ in terms of atoms, Studia Math. 62 (1978), 93–-101.

\bibitem{Stein}
    E. Stein,
    Harmonic analysis: real-variable methods, orthogonality, and oscillatory integrals,
    Princeton University Press, Princeton, NJ, 1993.

\bibitem{Stein-Weiss}
    E. Stein, G. Weiss,
    On the theory of harmonic functions of several variables. I. The theory of $H^p$-spaces.
    Acta Math. 103 (1960), 25–-62.
\end{thebibliography}
\end{document}